\documentclass[12pt]{amsart}
\usepackage{amssymb}
\usepackage{amscd}
\usepackage{amsthm}
\usepackage[dvips]{graphicx}
\usepackage[all,cmtip]{xy}
\usepackage{hyperref}

\usepackage{ulem}

\usepackage[all]{xy}

\newtheorem{Thm}{Theorem}[section]

\newtheorem{Def}[Thm]{Definition}

\newtheorem{Prop}[Thm]{Proposition}
\newtheorem{Ex1}[Thm]{Example}
\newtheorem{Rem1}[Thm]{Remark}

\newtheorem{Ques}[Thm]{Question}

\newenvironment{Rem}{\begin{Rem1}\rm}{\end{Rem1}}

\title{Two questions on stable equivalences of Morita type}
\author{Yuming Liu, Guodong Zhou and Alexander Zimmermann}

\address{Yuming Liu
\newline School of Mathematical Sciences
\newline Laboratory of Mathematics and Complex Systems
\newline Beijing Normal University
\newline Beijing 100875
\newline P.R.China}
\email{ymliu@bnu.edu.cn}

\address{Guodong Zhou
\newline Department of Mathematics
\newline Shanghai Key laboratory of PMMP
\newline East China Normal University
\newline  Dong Chuan Road 500
\newline Shanghai 200241
\newline P.R.China}
\email{gdzhou@math.ecnu.edu.cn}

\address{Alexander Zimmermann
\newline Universit\'e de Picardie,
\newline Facult\'e de Math\'ematiques et LAMFA (UMR 7352 du CNRS),
\newline 33 rue St Leu,
\newline F-80039 Amiens Cedex 1,
\newline France}
\email{Alexander.Zimmermann@u-picardie.fr}

  \date{version of \today}

\newenvironment{Proof}[1][Proof]{\begin{trivlist}
\item[\hskip \labelsep {\bfseries #1}]}{\flushright
$\Box$\end{trivlist}}

\newcommand{\lra}{\longrightarrow}

\newcommand{\ra}{\rightarrow}
\newcommand{\sdp}{\times\kern-.2em\vrule height1.1ex depth-.05ex}
\newcommand{\epi}{\lra \kern-.8em\ra}

\newcommand{\stmod}{\underline{\mathrm{mod}}}

 \normalbaselines
\thanks{The authors are supported by the exchange program STIC-Asie 'ESCAP' financed by the French Ministry of Foreign Affairs. The first author is supported by NCET Program from MOE of
China and by NNSF (No.11171325, No.11331006). The second author is supported    by Shanghai Pujiang
Program (No.13PJ1402800),  by National Natural Science Foundation of
China (No.11301186) and by the Doctoral Fund of Youth Scholars of
Ministry of Education of China (No.20130076120001).}

\begin{document}
\renewcommand{\thefootnote}{\alph{footnote}}
\renewcommand{\thefootnote}{\alph{footnote}}
\setcounter{footnote}{-1} 
\renewcommand{\thefootnote}{\alph{footnote}}
\setcounter{footnote}{-1} 

\begin{abstract}
It is well-known that derived equivalences preserve tensor products
and trivial extensions. We disprove both constructions for stable
equivalences of Morita type.
\end{abstract}

\maketitle

\section{Introduction}

\medskip
Let $k$ be a field and $A$ a finite dimensional $k$-algebra. We
denote by mod$A$ the category of all finite dimensional left
$A$-modules, and by $\stmod A$ the stable module category of mod$A$
modulo projective modules. Two finite dimensional $k$-algebras $A$
and $B$ are said to be {\it stably equivalent} if their stable
categories $\stmod A$ and $\stmod B$ are equivalent as
$k$-categories (\cite{ARS}). The stable category $\stmod A$ is a
natural quotient of the module category mod$A$ by the ideal of maps
that factor through projective modules, and in case that $A$ is
self-injective it is also a natural quotient (in the sense of
triangulated categories) of the bounded derived module category
$D^b$(mod$A$) (\cite{KeVo}\cite{Rickard1989b}). Examples of stable
equivalences naturally arise in the representation theory of groups
and algebras (see
\cite{ARS}\cite{AR1973}\cite{MV1980}\cite{Broue1994}\cite{Linckelmann1996}\cite{LiuXi2007}).

However, unlike the classical Morita theory for module categories
and the Morita theory for derived categories (\cite{Rickard1989a}),
it is not known how to describe stable equivalences in terms of
generators of stable categories (cf. \cite{KoenigLiu}). For this
reason, much less is known for stable equivalences comparing to
Morita and derived equivalences. In practice, one often uses {\it
stable equivalences of Morita type}, which form the class of stable
equivalences needed in most applications and which are close to
derived equivalences.

\begin{Def} \label{StM} $($\cite{Broue1994}$)$ Two finite dimensional
algebras $A$ and $B$ are said to be stably equi\-valent of Morita type if there are two
bimodules $_AM_B$ and $_BN_A$ which are projective as left modules
and as right modules such that there are bimodule isomorphisms:
$$ {}_A M\otimes_B N_A\simeq {}_AA_A\oplus {}_AP_A,\ \ \ {}_B
N\otimes_A M_B\simeq {}_BB_B\oplus {}_BQ_B$$ where ${}_AP_A$ and
${}_BQ_B$ are projective bimodules.
\end{Def}

Clearly, in the above situation, the exact functors $N\otimes_A-$
and $M\otimes_B-$ induce mutually inverse equivalences between $\stmod A$ and
$\stmod B$. In fact, any stable equivalence that is induced by an
exact functor between the module categories of two self-injective
algebras is isomorphic to a stable equivalence of Morita type
(\cite{Rickard1998}); Under some mild condition, this even holds for
general finite dimensional algebras (\cite{DugasMV}). All derived
equivalences between self-injective $k$-algebras induce stable
equivalences of Morita type (\cite{Rickard1991}). On the other hand,
there do exist stable equivalences of Morita type which are not
induced by derived equivalences (see
\cite{Broue1994}\cite{Linckelmann1996} and Section 3).

Although we have a better understanding of stable equivalences of
Morita type than general stable equivalences, we still cannot
answer some basic questions on them. For example, the fundamental
conjecture of Auslander and Reiten (see \cite{ARS}\cite{MV1990}),
which predicts that two stably equivalent algebras have the same
number of isomorphism classes of non-projective simple modules, is
largely open even for stable equivalences of Morita type (see
\cite{LiuZhouZimmermann}\cite{KoenigLiuZhou} for some equivalent
descriptions of this conjecture in this situation). In the present
paper, we will address another two basic questions on stable
equivalences of Morita type. Before stating these questions, we
first recall two classical results of Rickard on derived
equivalences.

\begin{Thm} $($\cite{Rickard1989b}\cite{Rickard1991}$)$ Let $A$ and $B$
be two derived equivalent finite dimensional $k$-algebras and assume the same condition for $C$ and $D$. Then
\begin{enumerate}
\item the trivial extension algebras $T(A)$ and $T(B)$ are derived equivalent;
\item the tensor product algebras $A\otimes_k C$ and  $B\otimes_k D$ are
derived equivalent.
\end{enumerate}
\end{Thm}

It is natural to ask whether the same are true for stable equivalences of Morita type. In fact,
such questions are closely related to the Auslander-Reiten conjecture.

\begin{Prop} \label{cyclicgroup} Let $k$ be an algebraically closed field of
characteristic $p>0$ and let $C_p$ be the cyclic group of order $p$.
Let $A$ and $B$  be two indecomposable, non-semisimple finite
dimensional algebras which are stably equivalent of Morita type.
Then the assertion that $A\otimes_k kC_p$ and  $B\otimes_k kC_p$ are
stably equivalent of Morita type implies the validity of the
Auslander-Reiten conjecture for $A$ and $B$.
\end{Prop}

\begin{Proof}(Compare with the proof of \cite[Theorem 3.7]{Rickard1998}) We first notice that
$A$ and $A\otimes_k kC_p$ have the same number of non-isomorphic simple modules.
Let $C_A$ be the Cartan matrix of $A$. The Cartan matrix of
$A\otimes_k kC_p$ is equal to $pC_A$, so its $p$-rank is zero. The
statement follows from Theorem 4.1 of \cite{LiuZhouZimmermann} which
says that the invariance of the $p$-rank of the Cartan matrix under
a stable equivalence of Morita type is equivalent to the
Auslander-Reiten conjecture.

One can also give a proof by computing the degree zero stable
Hochschild homology (see \cite{LiuZhouZimmermann}) of $A\otimes_k kC_p$ and of $B\otimes_k kC_p$.
The details are left to the reader.
\end{Proof}

In \cite{Rickard1998}, Rickard raised the following question.

 \begin{Ques} \label{Rickardquestion} $($\cite{Rickard1998}$)$ Let $A$ and $B$ be two indecomposable, non-semisimple self-injective $k$-algebras which are stably equivalent of Morita type and assume the same condition for $C$ and $D$. Are $A\otimes_kC$ and $B\otimes_kD$ stably equivalent of Morita type?
\end{Ques}

There would be trivial counterexamples if we do not request that algebras are indecomposable, since $A$ and $A\times k$ are stably equivalent of Morita type. If the stable equivalences are all induced by derived equivalences, then the answer is ``yes" since the derived equivalence preserves tensor product. If they are not, Rickard
mentions that the answer is probably ``no" in general. However, as Rickard stated, the simplest potential counterexamples are already quite complicated.

Note that in case $p=2$ and $A$ is symmetric, the above construction is just the
trivial extension, as the following more general proposition shows.

\begin{Prop} \label{trivialextension} Let $k$ be a field and $A$ be a symmetric $k$-algebra.
Then the tensor algebra $A\otimes_k k[x]/(x^2)$ is isomorphic to the trivial extension algebra $T(A)=A\ltimes D(A)$ of $A$.
\end{Prop}

\begin{Proof} Recall that as a $k$-vector space, $T(A)=A\oplus D(A)$, and the multiplication is defined by
$$(a,f)(b,g)=(ab, ag+fb)$$
for $a,b\in A$ and $f,g\in D(A)=\text{Hom}_k(A,k)$. Since $A$ is symmetric, we can fix an $A$-$A$-bimodule isomorphism
$A\rightarrow D(A)$ (write the image of $a$ by $a'$).
Define a map $$\alpha: A\otimes_k k[x]/(x^2)\rightarrow T(A)\mbox{ by }\alpha(a\otimes \overline{1} +b\otimes \overline{x})=(a,b').$$
It is straightforward that $\alpha$ is an algebra isomorphism.
\end{Proof}

\begin{Rem} \label{trivialextension-remark} Notice that one can define the trivial extension algebra $T(A)$
of arbitrary finite dimensional $k$-algebra $A$. It is well-known that $T(A)$ is always a symmetric $k$-algebra,
that is, $T(A)\cong D(T(A))$ as $T(A)$-$T(A)$-bimodules.
\end{Rem}

In \cite{KoenigLiuZhou}, K\"onig and the first named two authors
proved the following result relating the Auslander-Reiten conjecture
to trivial extensions.

\begin{Prop} \label{AR-trivialextensions} $($\cite[Corollary 8.2]{KoenigLiuZhou}$)$
Let $A$ and $B$ be two symmetric $k$-algebras over an algebraically
closed field of characteristic $p>0$. Suppose that $A$ and $B$ are
stably equivalent of Morita type. Then the condition that $T(A)$ and
$T(B)$ are stably equivalent of Morita type implies the validity of
the Auslander-Reiten conjecture for $A$ and $B$.
\end{Prop}

This motivates the following question in \cite{KoenigLiuZhou}.

\begin{Ques}  \label{KLZquestion} $($\cite[Question 8.3]{KoenigLiuZhou}$)$ Let $A$ and $B$ be two indecomposable, non-simple
finite dimensional algebras which are stably equivalent of Morita
type. Are their trivial extensions algebras $T(A)$ and
$T(B)$ stably equivalent of Morita type?
\end{Ques}

In the present paper, we will answer Question \ref{Rickardquestion}
and Question \ref{KLZquestion} to the negative for general finite
dimensional algebras. In particular, our result suggests that
Rickard's original question could have a negative answer in general.

\bigskip
\section{Triangular matrix algebras}
\medskip

In this section, we answer Question \ref{Rickardquestion} negatively
for general finite dimensional algebras.

Recall that for a finite dimensional $k$-algebra $A$, the stable
category $\overline{\text{mod}} A$ of mod$A$ modulo injective
modules can be defined similarly. There is an equivalence functor
$\tau$ from $\stmod A$ to $\overline{\text{mod}}A$, which is called
the Auslander-Reiten translation. If $F: \stmod A \ra \stmod B$ is a
stable equivalence, then there is an induced stable equivalence
 (modulo injectives) $\tau_B F \tau_A^{-1}: \overline{\text{mod}}A \ra
 \overline{\text{mod}}B$.

Given a finite dimensional $k$-algebra $A$, we denote by $T_2(A)$
the lower triangular matrix algebra $\left(
\begin{array}{cc} A&0\\ A&A \end{array}\right)$. Note that there is an algebra isomorphism between the tensor algebra $A\otimes_kT_2(k)$ and $T_2(A)$
under the map $a\otimes\left(
\begin{array}{cc} u&0\\ v&w \end{array}\right)\mapsto \left(
\begin{array}{cc} au&0\\ av&aw \end{array}\right)$. We refer to \cite{ARS} for the description of $T_2(A)$-modules in terms of $A$-modules.

\begin{Prop} \label{triangularmatrixalgebra}
Let $A$ and $B$ be two self-injective algebras with no semisimple
summands. If $\Lambda:= T_2(A)$ and
$\Gamma:= T_2(B)$
are stably equivalent, then $A$ and $B$ are Morita equivalent.
\end{Prop}

\begin{Proof} First we note that although $A$ and $B$ are self-injective algebras, $\Lambda$ and
$\Gamma$ are not self-injective any more. Suppose now that there is
a stable equivalence $F: \stmod \Lambda \ra \stmod \Gamma$. Let $H =
\tau_\Gamma F \tau_\Lambda^{-1}: \overline{\text{mod}}\Lambda \ra
\overline{\text{mod}}\Gamma$ be the induced
 stable equivalence modulo injectives. By \cite[Corollary
3.2]{AR1973}, $H$ induces a one-to-one correspondence between the
isomorphism classes of indecomposable non-simple non-injective
projective modules in mod$\Lambda$ and that in mod$\Gamma$. Under
our assumption, there are no simple projective modules over
$\Lambda$ and $\Gamma$. Therefore $H$ induces a one-to-one
correspondence between the isomorphism classes of indecomposable
non-injective projective modules in mod$\Lambda$ and that in
mod$\Gamma$.

Each $\Lambda$-module can be described as a triple $(X,Y,f)$, where
$X$ and $Y$ in mod$A$, and $f$ is an $A$-homomorphism from $X$ to
$Y$. A homomorphism from $(X,Y,f)$ to $(X',Y',f')$ is precisely a pair
$(\alpha, \beta)$ in $\text{Hom}_A(X, X') \times \text{Hom}_A(Y, Y')$ such that
$\beta f= f'\alpha.$ From this description we see that the
indecomposable projective $\Lambda$-modules are isomorphic to
modules of the form $(P,P,1_P)$ and $(0,P,0)$ where $P$ is an
indecomposable projective $A$-module. Dually, the indecomposable
injective $\Lambda$-modules are isomorphic to modules of the form
$(P,P,1_P)$ and $(P,0,0)$ where $P$ is an indecomposable projective
$A$-module. By the previous discussion, we know that under the
stable equivalence $H$, each indecomposable non-injective projective
$\Lambda$-module $(0,P,0)$ corresponds to some indecomposable
non-injective projective $\Gamma$-module $(0,Q,0)$, and this gives a
bijection between the isomorphism classes of indecomposable
non-injective projective modules in mod$\Lambda$ and that in
mod$\Gamma$. Observe that we have the following easy fact: for any
two $A$-modules $X$ and $X'$, we have
$\overline{\text{Hom}}_\Lambda((0,X,0),(0,X',0))\cong
\text{Hom}_\Lambda((0,X,0),(0,X',0))\cong\text{Hom}_A(X,X')$. Without loss of
generality we can assume that both $A$ and $B$ are basic algebras.
Then we have that $$H((0,A,0))\cong (0,B,0) \mbox{ and
}\overline{\text{End}}_\Lambda((0,A,0))\cong
\overline{\text{End}}_\Gamma((0,B,0)).$$ Therefore we have the following
algebra isomorphisms:
$$\text{End}_A(A)\cong \text{End}_\Lambda((0,A,0))\cong
\overline{\text{End}}_\Lambda((0,A,0))\cong \quad \quad \quad \quad \quad
\quad  \quad \quad  \quad \quad $$
$$\quad \quad \quad \quad \quad \quad  \quad \quad  \quad \quad \cong \overline{\text{End}}_\Gamma((0,B,0))\cong \text{End}_\Gamma((0,B,0))\cong \text{End}_B(B).$$
It follows that the algebras $A$ and $B$ are isomorphic.
\end{Proof}

\begin{Rem} The above result shows that the answer to Question \ref{Rickardquestion} is negative for general finite dimensional algebras. Indeed, we can easily
find two self-injective algebras $A$ and $B$ which are derived
equivalent but not Morita equivalent. Clearly $A$ and $B$ are stably
equivalent of Morita type, but $T_2(A)\simeq A\otimes_kT_2(k)$ and $T_2(B)\simeq B\otimes_kT_2(k)$
cannot be stably equivalent of Morita type by Proposition \ref{triangularmatrixalgebra}.
\end{Rem}

\begin{Rem} From the proof of Proposition \ref{triangularmatrixalgebra}, we obtain that the stable category of the triangular matrix algebra $T_2(A)$ determines the
original algebra $A$ in the following way: it is the (stable) endomorphism algebra of the sum of indecomposable non-projective injective
modules over triangular matrix algebra.
\end{Rem}

\bigskip
\section{Trivial extensions}
\medskip

 In this section, we answer Question \ref{KLZquestion} negatively for general finite dimensional algebras.

Let $k$ be an algebraically closed field of characteristic $2$. Then it is well-known (see for example, \cite{Erdmann}) that the group algebra $A=kD_8$ of the dihedral group of order $8$ is given by the following quiver
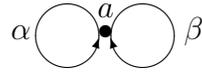
\unitlength=1.00mm \special{em:linewidth 0.4pt}
\linethickness{0.4pt}
\begin{center}
\begin{picture}(32.00,12.00)
\put(18,5){$\bullet$}\put(14,5.5){\circle{8}}\put(17.7,3.8){\vector(1,3){0.5}}\put(6.5,5){$\alpha$}
\put(24,5.5){\circle{8}}\put(20.2,3.8){\vector(-1,3){0.5}}\put(29.5,5){$\beta$}\put(18,8){$a$}
\end{picture}
\end{center}
with relations
$$\qquad \qquad \alpha^2= \beta^2 =0, (\alpha\beta)^2=(\beta\alpha)^2.$$
This is a local symmetric algebra with basis (for simplicity, we write $\alpha$ for $\overline{\alpha}$, and
 so on) $1=e_a,\alpha,\beta,\alpha\beta,\beta\alpha,\alpha\beta\alpha,\beta\alpha\beta,
\alpha\beta\alpha\beta=\beta\alpha\beta\alpha$. The Loewy diagram of the regular module $_AA$ looks like
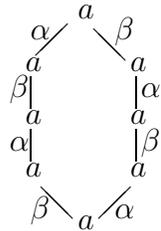
\unitlength=0.7cm \special{em:linewidth 0.4pt}
\linethickness{0.4pt}
\begin{center}
\begin{picture}(6.00,6.00)
\put(1,3){$a$}
\put(3,3){$a$} \put(1,4){$a$}\put(3,4){$a$}
\put(2,5){$a$} \put(1.2,4.3){\line(1,1){0.6 }}
\put(2.3,4.8){\line(1,-1){0.6}}\put(1.1,4.6){$\alpha$}\put(2.7,4.6){$\beta$}
\put(1.1,3.9){\line(0,-1){0.6}}\put(3.1,3.9){\line(0,-1){0.6}}\put(0.7,3.5){$\beta$}\put(3.2,3.5){$\alpha$}
\put(1.1,2.9){\line(0,-1){0.6}}\put(3.1,2.9){\line(0,-1){0.6}}\put(0.7,2.5){$\alpha$}\put(3.2,2.5){$\beta$}
\put(1,2){$a$}\put(3,2){$a$}
\put(2,1){$a$}
 \put(1.3,1.8){\line(1,-1){0.6 }}
\put(2.4,1.2){\line(1,1){0.6}}\put(1.1,1.2){$\beta$}\put(2.7,1.2){$\alpha$}
\end{picture}
\end{center}
The Cartan matrix of $A$ is given by $C_A=(8)$, and the center
$Z(A)$ of $A$ is a radical square zero local algebra with basis
$1,\alpha\beta+\beta\alpha,\alpha\beta\alpha,\beta\alpha\beta,\alpha\beta\alpha\beta$. Let
$S$ be the unique simple $A$-module (which is also the trivial
module $k$ of the group algebra $A$) and let $\Lambda$ be the
endomorphism algebra End$_A(A\oplus S)^{op}$. One can compute that
$\Lambda$ is given by the following quiver \unitlength=1.00mm
\special{em:linewidth 0.4pt} \linethickness{0.4pt}
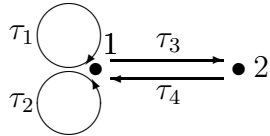
\begin{center}
\begin{picture}(50.00,16.00)
\put(18,5){$\bullet$}\put(15.5,10.5){\circle{8}}\put(18.2,7.2){\vector(-3,-4){0.5}}\put(7.5,10){$\tau_1$}
\put(15.5,1.5){\circle{8}}\put(19,3.8){\vector(-1,2){0.5}}\put(7.5,1){$\tau_2$}\put(20,8){$1$}
\put(21,7.2){\vector(1,0){15}}\put(20,8){$1$}\put(27,9){$\tau_3$}
\put(36,4.8){\vector(-1,0){15}}\put(27,2){$\tau_4$}
\put(37,5){$\bullet$}\put(40,5){$2$}
\end{picture}
\end{center}
with relations
$${\tau_1}^2= {\tau_2}^2 =\tau_3\tau_4=\tau_2\tau_4=\tau_1\tau_4=\tau_3\tau_1=\tau_3\tau_2=0, (\tau_1\tau_2)^2=(\tau_2\tau_1)^2=\tau_4\tau_3.$$
This is an $11$-dimensional algebra with basis (remind that we write $\tau_1$ for $\overline{\tau_1}$, and so on) $$e_1,e_2,\tau_1,\tau_2,\tau_3,\tau_4,\tau_2\tau_1,\tau_1\tau_2,\tau_1\tau_2\tau_1,
\tau_2\tau_1\tau_2,\tau_2\tau_1\tau_2\tau_1=\tau_1\tau_2\tau_1\tau_2=\tau_4\tau_3. $$
The regular module $_\Lambda\Lambda$ has the following decomposition
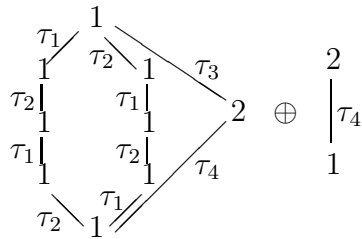
\unitlength=0.7cm \special{em:linewidth 0.4pt}
\linethickness{0.4pt}
\begin{center}
\begin{picture}(12.00,6.00)
\put(1,3){$1$}
\put(3,3){$1$} \put(1,4){$1$}\put(3,4){$1$}
\put(2,5){$1$} \put(2.5,5){\line(3,-2){2.1}}\put(4.6,3.2){\line(-1,-1){2.1}}\put(4.7,3.2){$2$}
\put(4.0,4.1){$\tau_3$}\put(4.0,2.2){$\tau_4$}
\put(1.2,4.3){\line(1,1){0.6 }}
\put(2.3,4.8){\line(1,-1){0.6}}\put(1.0,4.7){$\tau_1$}\put(2,4.3){$\tau_2$}
\put(1.1,3.9){\line(0,-1){0.5}}\put(3.1,3.9){\line(0,-1){0.5}}\put(0.5,3.5){$\tau_2$}\put(2.5,3.5){$\tau_1$}
\put(1.1,2.9){\line(0,-1){0.5}}\put(3.1,2.9){\line(0,-1){0.5}}\put(0.5,2.5){$\tau_1$}\put(2.5,2.5){$\tau_2$}
\put(1,2){$1$}\put(3,2){$1$}
\put(2,1){$1$}
 \put(1.3,1.8){\line(1,-1){0.6 }}
\put(2.4,1.2){\line(1,1){0.6}}\put(1.0,1.2){$\tau_2$}\put(2.2,1.6){$\tau_1$}
\put(5.5,3.2){$\oplus$}\put(6.5,4.2){$2$}\put(6.5,2.2){$1$}
\put(6.6,4.0){\line(0,-1){1.2}}\put(6.7,3.2){$\tau_4$}
\end{picture}
\end{center}
The Cartan matrix of $\Lambda$ is given by $C_\Lambda=\left(
\begin{array}{cc} 8&1\\ 1&1 \end{array}\right)$, and the center $Z(\Lambda)$ of $\Lambda$ is a $5$-dimensional algebra with basis $1,\tau_2\tau_1+\tau_1\tau_2,\tau_1\tau_2\tau_1,\tau_2\tau_1\tau_2,\tau_2\tau_1\tau_2\tau_1=\tau_1\tau_2\tau_1\tau_2=\tau_4\tau_3$.
Since char$k=2$, it is easy to verify that $Z(\Lambda)$ is also a radical square zero local algebra.
Next we want to compute the center $Z(T(\Lambda))$ of the trivial extension $T(\Lambda)=\Lambda\ltimes D(\Lambda)$.
According to a result of Bessenrodt, Holm and the third named author (see \cite[Proposition 3.2]{BHZ}),
$Z(T(\Lambda))=Z(\Lambda)\ltimes \text{Ann}_{{D(\Lambda)}}(K(\Lambda))$, where $K(\Lambda)$ is the
$k$-subspace of $\Lambda$ spanned by all commutators $\lambda\mu - \mu\lambda$ $(\lambda,\mu\in \Lambda)$
and where $\text{Ann}_{{D(\Lambda)}}(K(\Lambda))=\{f\in {D(\Lambda)}|f(K(\Lambda))=0\}$. By a straightforward
calculation we have the following (write ${e_1}^*$ as the dual basis element corresponding to $e_1$, and so on):
$$K(\Lambda)=\langle \tau_2\tau_1+\tau_1\tau_2,\tau_3,\tau_4,\tau_1\tau_2\tau_1,\tau_2\tau_1\tau_2,\tau_4\tau_3\rangle,$$
$$\text{Ann}_{{D(\Lambda)}}(K(\Lambda))=\langle {e_1}^*,{e_2}^*,{\tau_1}^*,{\tau_2}^*,(\tau_2\tau_1)^*+(\tau_1\tau_2)^*\rangle.$$
Again char $k=2$ forces that $Z(T(\Lambda))$ is a ($10$-dimensional)
radical square zero local algebra.

Now we come back to the group algebra $A$ of the dihedral group $D_8$. According to \cite{CT}, each direct summand of the following module
\unitlength=0.7cm \special{em:linewidth 0.4pt}
\linethickness{0.4pt}
\begin{center}
\begin{picture}(8.00,2.50)\put(-4,1){$\text{rad}A/\text{rad}^4A$} \put(-0.3,1){$=$}  \put(1.9,1){$\oplus$}
\put(1,1){$a$}
\put(3,1){$a$} \put(1,2){$a$}\put(3,2){$a$}\put(1,0){$a$}
\put(3,0){$a$} \put(1.1,1.9){\line(0,-1){0.6}}\put(3.1,1.9){\line(0,-1){0.6}}\put(0.6,1.5){$\beta$}\put(3.2,1.5){$\alpha$}
\put(1.1,0.9){\line(0,-1){0.6}}\put(3.1,0.9){\line(0,-1){0.6}}\put(0.6,0.5){$\alpha$}\put(3.2,0.5){$\beta$}
\end{picture}
\end{center}
is an endotrivial module. Recall that for a group algebra $kG$ of a
finite group $G$, a $kG$-module $X$ is called endotrivial if
$D(X)\otimes_kX\simeq k\oplus \mbox{\{projective\}}$ as $kG$-modules
(where $k$ is the trivial module). It follows easily that
$X\otimes_k-$ induces a stable self-equivalence of Morita type over
$\stmod kG$ (here the defining $kG$-$kG$-bimodule is given by
$X\otimes_kkG$ where the left $kG$-module structure is defined by
diagonal $G$-action and the right $kG$-module structure is defined
by by multiplication on the right factor). Let now $X$ be one of any
direct summand of $\text{rad}A/\text{rad}^4A$. Then $X$ induces a
stable self-equivalence of Morita type over $\stmod A$ such that the
trivial module $S$ corresponds to $X$. Since $X$ is not of the form
$\Omega^i(S)$, this stable self-equivalence is not induced from a
derived self-equivalence (cf. \cite[Remark 3.10]{Linckelmann1996}
or \cite[Theorem 2.11]{TrPic}). Let $\Gamma$ be the endomorphism
algebra End$_A(A\oplus X)^{op}$. Then by the construction in
\cite[Theorem 1.1]{LiuXi2007}, there is a stable equivalence of
Morita type between $\Lambda$ and $\Gamma$. One can compute that
$\Gamma$ has the same quiver as $\Lambda$ (but here we use new
notations to name the arrows) \unitlength=1.00mm
\special{em:linewidth 0.4pt} \linethickness{0.4pt}
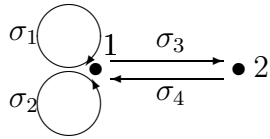
\begin{center}
\begin{picture}(50.00,16.00)
\put(18,5){$\bullet$}\put(15.5,10.5){\circle{8}}\put(18.2,7.2){\vector(-3,-4){0.5}}\put(7.5,10){$\sigma_1$}
\put(15.5,1.5){\circle{8}}\put(19,3.8){\vector(-1,2){0.5}}\put(7.5,1){$\sigma_2$}\put(20,8){$1$}
\put(21,7.2){\vector(1,0){15}}\put(20,8){$1$}\put(27,9){$\sigma_3$}
\put(36,4.8){\vector(-1,0){15}}\put(27,2){$\sigma_4$}
\put(37,5){$\bullet$}\put(40,5){$2$}
\end{picture}
\end{center}
with relations
$${\sigma_1}^2= {\sigma_2}^2 =\sigma_3\sigma_1=\sigma_2\sigma_4=0,\sigma_2\sigma_1=\sigma_4\sigma_3,$$
$$(\sigma_2\sigma_1)^2=(\sigma_1\sigma_2)^2=(\sigma_4\sigma_3)^2=\sigma_1\sigma_4\sigma_3\sigma_2,
\sigma_3\sigma_2\sigma_1\sigma_2=0.$$
This is a $16$-dimensional algebra with basis $$e_1,e_2,\sigma_1,\sigma_2,\sigma_3,\sigma_4,\sigma_2\sigma_1=\sigma_4\sigma_3,
\sigma_1\sigma_2,\sigma_3\sigma_2,\sigma_3\sigma_4,\sigma_1\sigma_4,
\sigma_1\sigma_2\sigma_1,\sigma_2\sigma_1\sigma_2,$$
$$\sigma_3\sigma_2\sigma_1=\sigma_3\sigma_4\sigma_3,\sigma_4\sigma_3\sigma_4=\sigma_2\sigma_1\sigma_4,
(\sigma_2\sigma_1)^2=(\sigma_1\sigma_2)^2=(\sigma_4\sigma_3)^2. $$
The regular module $_\Gamma\Gamma$ has the following decomposition
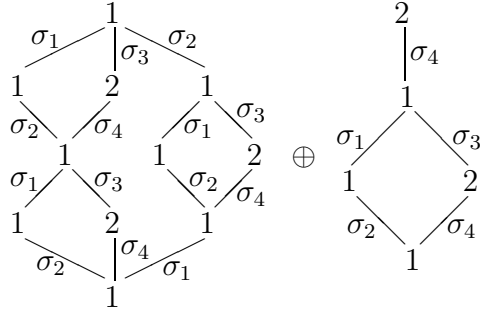
\unitlength=0.7cm \special{em:linewidth 0.4pt}
\linethickness{0.4pt}
\begin{center}
\begin{picture}(12.00,6.00)
\put(3,5){$1$}
 \put(3,4.9){\line(-2,-1){1.5}}\put(3.2,4.9){\line(0,-1){0.8}}\put(3.4,4.9){\line(2,-1){1.5}}
\put(1.6,4.6){$\sigma_1$}\put(3.3,4.3){$\sigma_3$}\put(4.2,4.6){$\sigma_2$}
\put(1.2,3.6){$1$}\put(3,3.6){$2$}\put(4.8,3.6){$1$}
 \put(1.4,3.5){\line(1,-1){0.7}}\put(3.1,3.5){\line(-1,-1){0.7}}\put(4.8,3.5){\line(-1,-1){0.7}}
 \put(5.1,3.5){\line(1,-1){0.7}}
\put(1.2,2.9){$\sigma_2$}\put(2.8,2.9){$\sigma_4$}\put(4.5,2.9){$\sigma_1$}\put(5.5,3.2){$\sigma_3$}
\put(2.1,2.3){$1$}\put(3.9,2.3){$1$}\put(5.7,2.3){$2$}
 \put(2.2,2.2){\line(-1,-1){0.7}}\put(2.4,2.2){\line(1,-1){0.7}}\put(4.2,2.2){\line(1,-1){0.7}}
 \put(5.8,2.2){\line(-1,-1){0.7}}
\put(1.2,1.9){$\sigma_1$}\put(2.8,1.9){$\sigma_3$}\put(4.6,1.9){$\sigma_2$}\put(5.5,1.6){$\sigma_4$}
\put(1.2,1){$1$}\put(3,1){$2$}\put(4.8,1){$1$}
 \put(1.5,0.9){\line(2,-1){1.5}}\put(3.2,0.9){\line(0,-1){0.8}}\put(4.9,0.9){\line(-2,-1){1.5}}
\put(1.7,0.3){$\sigma_2$}\put(3.3,0.6){$\sigma_4$}\put(4.1,0.2){$\sigma_1$}
\put(3,-0.4){$1$}
\put(6.55,2.3){$\oplus$}
\put(8.5,5){$2$}
\put(8.7,4.9){\line(0,-1){1.0}}
\put(8.8,4.3){$\sigma_4$}
\put(8.6,3.4){$1$}
\put(8.7,3.3){\line(-1,-1){1.0}}\put(8.9,3.3){\line(1,-1){1.0}}
\put(7.4,2.7){$\sigma_1$}\put(9.6,2.7){$\sigma_3$}
\put(7.5,1.8){$1$}\put(9.8,1.8){$2$}
\put(7.8,1.7){\line(1,-1){0.9}}\put(9.9,1.7){\line(-1,-1){0.9}}
\put(7.6,1.0){$\sigma_2$}\put(9.5,1.0){$\sigma_4$}
\put(8.7,0.3){$1$}
\end{picture}
\end{center}
\vspace{0.5cm}
The Cartan matrix of $\Gamma$ is given by $C_\Gamma=\left(
\begin{array}{cc} 8&3\\ 3&2 \end{array}\right)$. From this we can deduce that $\Lambda$ and $\Gamma$ are not derived equivalent,
since their Cartan matrices are not congruent over the integers. The fact that the Cartan matrices of two derived equivalent algebras
are congruent over the integers is one of few known invariants to distinguish between derived equivalence and stable equivalence of Morita type.
Our next aim is to show that the trivial extensions $T(\Lambda)$ and $T(\Gamma)$ are also not stably equivalent of Morita type.
We will verify this fact by proving that their stable centers are not isomorphic as algebras.

Let us first recall the definition of the stable center. For an algebra $A$, we can identify any $A$-$A$-bimodule with a left
$A^e$-module where $A^e=A\otimes_kA^{op}$. In particular, the algebra $A$ itself is naturally an $A^e$-module, and the endomorphism
algebra $\mathrm{End}_{A^e}(A,A)$ is canonically isomorphic to the center $Z(A)$ of $A$ $(\mbox{by } f\mapsto f(1))$.
Set $Z^{{pr}}(A)$ to be the ideal of $Z(A)$ consisting of
$A^e$-homomorphisms from $A$ to $A$ which factor through a
projective $A^e$-module and we call it the {\it projective center} of $A$.
The {\it stable center} of $A$ is defined to be the quotient algebra $Z^{st}(A)=Z(A)/Z^{{pr}}(A)$. It is well-known that a stable
equivalence of Morita type preserves the stable centers of algebras (see \cite{Broue1994}).

\begin{Prop} \label{nonisomorphicstablecenters}
Let $\Lambda$ and $\Gamma$ be as above. Then the stable centers of $T(\Lambda)$ and $T(\Gamma)$ are not isomorphic as algebras.
In particular, $T(\Lambda)$ and $T(\Gamma)$ are not stably equivalent of Morita type.
\end{Prop}

\begin{Proof} The Cartan matrices of $\Lambda$ and $\Gamma$ are given by $C_\Lambda=\left(
\begin{array}{cc} 8&1\\ 1&1 \end{array}\right)$ and $C_\Gamma=\left(
\begin{array}{cc} 8&3\\ 3&2 \end{array}\right)$, respectively. It follows easily that the Cartan matrices of $T(\Lambda)$
and $T(\Gamma)$ are given by $C_{T(\Lambda)}=\left(
\begin{array}{cc} 16&2\\ 2&2 \end{array}\right)$ and $C_{T(\Gamma)}=\left(
\begin{array}{cc} 16&6\\ 6&4 \end{array}\right)$, respectively. By \cite[Proposition 2.3 and Corollary 2.7]{LiuZhouZimmermann},
the dimension of the projective center of a symmetric algebra over an algebraically closed field $k$ of characteristic $p\geq 0$
is equal to the $p$-rank of the Cartan matrix. Since now $p=2$, both $2$-ranks of $C_{T(\Lambda)}$ and $C_{T(\Gamma)}$ are zero,
and therefore the stable centers of $T(\Lambda)$ and $T(\Gamma)$ are the same as the centers of $T(\Lambda)$ and $T(\Gamma)$, respectively.

We have seen that the center $Z(T(\Lambda))$ is a $10$-dimensional radical square zero local algebra.
Similarly we can compute the center $Z(T(\Gamma))$ using the formula $Z(T(\Gamma))=Z(\Gamma)\ltimes \text{Ann}_{{D(\Gamma)}}(K(\Gamma))$.
The center $Z(\Gamma)$ of $\Gamma$ is a $5$-dimensional algebra with basis
$1,\sigma_2\sigma_1+\sigma_1\sigma_2+\sigma_3\sigma_4,\sigma_1\sigma_2\sigma_1,
\sigma_2\sigma_1\sigma_2,(\sigma_2\sigma_1)^2$.
Since char$k=2$, it is easy to verify that $Z(\Gamma)$ is also a radical square zero local algebra. We also have the following:
$$K(\Gamma)=\langle \sigma_3,\sigma_4,\sigma_1\sigma_4,\sigma_3\sigma_2,
\sigma_3\sigma_2\sigma_1=\sigma_3\sigma_4\sigma_3,\sigma_4\sigma_3\sigma_4,
\sigma_2\sigma_1+\sigma_1\sigma_2,\sigma_3\sigma_4+\sigma_4\sigma_3,$$
$$\sigma_2\sigma_1\sigma_2,\sigma_1\sigma_2\sigma_1=\sigma_1\sigma_4\sigma_3,
(\sigma_2\sigma_1)^2=(\sigma_1\sigma_2)^2=(\sigma_4\sigma_3)^2\rangle,$$
$$\text{Ann}_{{D(\Gamma)}}(K(\Gamma))=\langle {e_1}^*,
{e_2}^*,{\sigma_1}^*,{\sigma_2}^*,(\sigma_2\sigma_1)^*+(\sigma_1\sigma_2)^*+(\sigma_3\sigma_4)^*\rangle.$$
We perform the following multiplication in $Z(T(\Gamma))$:
$$(\sigma_2\sigma_1+\sigma_1\sigma_2+\sigma_3\sigma_4)((\sigma_2\sigma_1)^*+(\sigma_1\sigma_2)^*+(\sigma_3\sigma_4)^*)=2{e_1}^*+{e_2}^*={e_2}^*.$$
Since char $k=2$, the above multiplication is not equal to zero and
therefore $Z(T(\Gamma))$ is not radical square zero. So
$Z(T(\Lambda))$ and $Z(T(\Gamma))$ are not isomorphic as algebras.

\end{Proof}

\begin{Rem} Suppose that $k$ is of characteristic $2$.
Then the center $Z(T(\Gamma))$ is a $10$-dimensional local algebra such that the regular module has the following Loewy structure
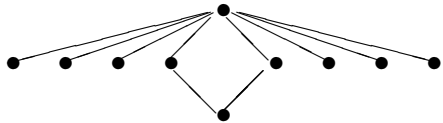
\unitlength=0.7cm \special{em:linewidth 0.4pt}
\linethickness{0.4pt}
\begin{center}
\begin{picture}(12.00,3.00)
\put(5,0){$\bullet$}
\put(1,1){$\bullet$}\put(2,1){$\bullet$}\put(3,1){$\bullet$}\put(4,1){$\bullet$}\put(6,1){$\bullet$}
\put(7,1){$\bullet$}\put(8,1){$\bullet$}\put(9,1){$\bullet$}
\put(5,2){$\bullet$}
\put(4.9,2){\line(-1,-1){.8}}
\put(5.2,2){\line(1,-1){.8}}
\put(4.2,1){\line(1,-1){.8}}
\put(5.9,1){\line(-1,-1){.8}}
\put(5.9,1){\line(-1,-1){.8}}
\put(5.3,2.1){\line(2,-1){1.8}}
\put(5.4,2.1){\line(3,-1){2.8}}
\put(5.45,2.1){\line(4,-1){3.8}}
\put(4.95,2.1){\line(-2,-1){1.8}}
\put(4.9,2.1){\line(-3,-1){2.8}}
\put(4.85,2.1){\line(-4,-1){3.8}}
\end{picture}
\end{center}
Parallel edges correspond to multiplication with the same element.
Here, in the square in the centre one direction
corresponds to multiplication with
$(\sigma_2\sigma_1+\sigma_1\sigma_2+\sigma_3\sigma_4)$, whereas the
other direction of the square in the centre corresponds to
multiplication with
$(\sigma_2\sigma_1)^*+(\sigma_1\sigma_2)^*+(\sigma_3\sigma_4)^*$.
The product of these two corresponds to multiplication with $e_2^*$.
\end{Rem}

\end{document}